\documentclass[10pt]{amsart}

\usepackage{amstext,amsfonts,amssymb,amscd,amsbsy,amsmath}
\usepackage{ifthen}
\usepackage{amsthm}
\usepackage{latexsym}
\usepackage[all]{xy}
\usepackage{enumerate}
\usepackage[mathscr]{eucal}

\theoremstyle{plain}
\newtheorem{theorem}{Theorem}[section]

\newtheorem{corollary}[theorem]{Corollary}
\newtheorem{proposition}[theorem]{Proposition}

\theoremstyle{plain}
\newtheorem*{question}{Open Question}

\theoremstyle{definition}
\newtheorem{definition}[theorem]{Definition}
\newtheorem{example}[theorem]{Example}

\theoremstyle{remark}
\newtheorem*{remark}{Remark}
\newtheorem*{acknowledge}{Acknowledgements}

\DeclareMathOperator{\Span}{Span}
\DeclareMathOperator{\Ad}{Ad}
\DeclareMathOperator{\red}{red}
\DeclareMathOperator{\dom}{dom}
\DeclareMathOperator{\ran}{ran}

\newcommand{\bh}{\mcB(\mcH)}
\newcommand{\el}[1]{\ell^{2}(#1)}
\newcommand{\eps}{\varepsilon}
\newcommand{\inp}[1]{\left<#1\right>}
\newcommand{\motimes}{{\otimes_{\mbox{\tiny{min}}}}}
\newcommand{\norm}[1]{\left\|#1\right\|}
\newcommand{\ua}{\uparrow\negmedspace\negmedspace}

\newcommand{\mcT}{\mathcal{T}}
\newcommand{\mcH}{\mathcal{H}}
\newcommand{\mcB}{\mathcal{B}}
\newcommand{\mcS}{\mathcal{S}}
\newcommand{\mcW}{\mathcal{W}}
\newcommand{\bbC}{\mathbb{C}}
\newcommand{\bbF}{\mathbb{F}}
\newcommand{\bbN}{\mathbb{N}}
\newcommand{\bbZ}{\mathbb{Z}}
\newcommand{\bbT}{\mathbb{T}}
\newcommand{\ra}{\rightarrow}

\newcommand{\into}{\hookrightarrow}
\newcommand{\onto}{\twoheadrightarrow}
\newcommand{\cstar}{\text{C}\sp{*}}
\newcommand{\msE}{\mathscr{E}}

\begin{document}

\title[amenability and weak containment]{$C^*$-algebras of inverse semigroups: amenability and weak containment}

\author{David Milan}
\address{David Milan, Department of Mathematics, University of Nebraska, {Lincoln, NE}
68588-0130, USA}
\email{dmilan@unl.math.edu}
\begin{abstract}
We argue that weak containment is an appropriate notion of am-enability for inverse semigroups. Given an inverse semigroup $S$ and a homomorphism $\varphi$ of $S$ onto a group $G$, we show, under an assumption on $\ker(\varphi)$, that $S$ has weak containment if and only if $G$ is amenable and $\ker(\varphi)$ has weak containment. Using Fell bundle amenability, we find a related result for inverse semigroups with zero. We show that all graph inverse semigroups have weak containment and that Nica's inverse semigroup $\mcT_{G,P}$ of a quasi-lattice ordered group $(G,P)$ has weak containment if and only if $(G,P)$ is amenable.
\end{abstract}
\maketitle

\section*{Introduction} \label{secIntro}
Amenability for inverse semigroups has been studied by a number of authors, and the results suggest that not all of the equivalent definitions of group amenability translate well to inverse semigroups. In \cite{DunNam78} it is shown that an inverse semigroup
$S$ admits a left invariant mean if and only if the maximum group homomorphic image of $S$ is amenable. This notion of amenability is too weak for inverse semigroups. Indeed, any inverse semigroup with zero has trivial maximum group homomorphic image and hence admits an invariant mean (the unique invariant mean is given by evaluating functions in $\ell^{\infty}(S)$ at zero). Another equivalent condition from group theory says that a group $G$ is amenable if and only if $L^{1}(G)$ is amenable as a Banach algebra. It is also shown in \cite{DunNam78} for a discrete inverse semigroup that $\ell^1(S)$ is amenable as a Banach algebra if and only if the set of idempotents of $S$ is finite and every subgroup of $S$ is amenable. This notion is too strong, for example, since many infinite commutative inverse semigroups would fail to be amenable.

As in the work of Duncan and Paterson \cite{DunPat85}, we study yet another notion of amenability motivated by group theory: the weak containment property. Recall that the full $C^*$-algebra $C^*(S)$ is universal for representations of $S$ by partial isometries. That is, any representation $\pi : S \ra PI(\mcH)$, where $PI(\mcH)$ is the set of partial isometries on a Hilbert space $\mcH$, induces a $*$-homomorphism $\pi : C^*(S) \ra \bh$. In particular, the left regular representation induces a $*$-homomorphism $\Lambda : C^*(S) \ra \mcB(\el{S})$ whose image is denoted $C^*_r(S)$ and is called the \textit{reduced $C^*$-algebra} of $S$. An inverse semigroup $S$ has \textit{weak containment} if and only if $\Lambda$ is an isomorphism.

In many ways, the weak containment property is an appropriate notion of amenability for inverse semigroups. For example, if $F_2$ denotes the free group on two generators, then $F_2^0$ admits an invariant mean (we denote by $S^0$ the inverse semigroup $S$ with an adjoined zero and the obvious multiplication). However, $F_2^0$ fails to have weak containment. On the other hand, any commutative inverse semigroup has weak containment.

We prove some additional results that suggest weak containment is the right notion of amenability for inverse semigroups. Let $\varphi: S \ra G$ be a homomorphism onto a group with kernel $H$. We say $H$ is \textit{$C^*$-isometric} in $S$ if $\norm{g}_{C^*(S)} = \norm{g}_{C^*H}$ for all $g$ in $\bbC H$. Assuming $H$ is $C^*$-isometric in $S$, we show $S$ has weak containment if and only if $G$ is amenable and $H$ has weak containment (Corollary \ref{cor_weakcontainment}). This resembles the theorem that says, for an exact sequence of discrete groups $H \into G \onto K$, that $G$ is amenable if and only if $H$ and $K$ are amenable.

Every inverse semigroup $S$ has an associated maximum group homomorphic image $G(S)$, where the homomorphism $\sigma : S \ra G(S)$ is the quotient map obtained by identifying $s,t$ in $S$ if $es = et$ for some idempotent $e$. We say $S$ is \textit{$E$-unitary} if $\ker \sigma$ consists solely of idempotents (in which case we say $\sigma$ is \textit{idempotent pure}). As a consequence of Corollary \ref{cor_weakcontainment}, an $E$-unitary inverse semigroup $S$ has weak containment if and only if $G(S)$ is amenable. This generalizes what is already known for $E$-unitary Clifford semigroups \cite{Pat78}.

Many inverse semigroups suffer from a seemingly fatal flaw: they contain a zero. In that case any group homomorphic image is the trivial group, and the results described above are vacuous. This difficulty is overcome in section 3 by replacing homomorphisms with
maps $\varphi : S \ra G^0$ such that $\varphi(ab) = \varphi(a) \varphi(b)$ whenever $ab \neq 0$ and $\varphi^{-1}(0) = \{0\}$. A map $\varphi$ satisfying these properties is called a \textit{grading} of $S$ by the group $G$. This approach is inspired by Lawson's construction of a grading $\sigma: S \ra G^{0}$, called the universal grading of $S$, that generalizes the maximum group homomorphism \cite{Lawson01}. Building on the work of
Bulman-Fleming, Fountain, and Gould \cite{BFG99}, Lawson shows that the so-called strongly $E^*$-unitary inverse semigroups are exactly the inverse semigroups having idempotent pure universal gradings. Thus the property of being strongly $E^*$-unitary is the right
generalization of the $E$-unitary property to inverse semigroups with zero.

Let $\varphi$ be a grading of an inverse semigroup $S$ by a group $G$ with kernel $H$. Unlike before, we can not say that $G$ is amenable if $S$ has weak containment. Instead, we relate the weak containment property of $S$ to the amenability of a $C^*$-algebraic bundle over $G$ induced from $\varphi$. (Such amenability was defined and studied by Exel \cite{Exel97}, where the term Fell bundle is used instead of $C^*$-algebraic bundle. We outline the results of Exel that will be needed in this paper in the next section.) As a consequence, a strongly $E^*$-unitary inverse semigroup $S$ has weak containment if and only if the associated Fell bundle over the universal group is amenable.

Using a result from \cite{Exel00} on Fell bundle amenability, we show that graph inverse semigroups have weak containment. It is known that graph inverse semigroups are strongly $E^*$-unitary with free universal gradings \cite{Lawson01}. Thus, for example, the polycyclic inverse semigroup on $n$ generators (the graph inverse semigroup associated with the bouquet of $n$ circles) has weak containment, even though its universal group is $F_n$.

We also consider the Toeplitz inverse semigroup $\mcT_{G,P}$ defined by Nica \cite{Nica94}. We note in section 5 that $\mcT_{G,P}$ has weak containment if and only if the pair $(G,P)$ is amenable in the sense defined by Nica \cite{Nica92}.

In the final section we study positivity of the restriction map from the complex algebra $\bbC S$ of an inverse semigroup onto the complex algebra $\bbC H$ of an inverse subsemigroup. This property is important because it implies that $H$ is $C^*$-isometric in $S$. We find some classes of inverse semigroups where the restriction map is positive, but we do not know if positivity holds whenever $H$ is the kernel of a homomorphism onto a group (or when $H$ is the kernel of a group grading).

\section{Preliminaries} \label{secBack}
A semigroup $S$ is an \textit{inverse semigroup} if for each $s$ in $S$ there exists unique $s^*$ in $S$ such that $s = s s^* s$ and $s^* = s^* s s^*$. In this paper we deal only with discrete inverse semigroups (and discrete groups). There is a natural partial order on $S$ defined by $s \leq t$ if $s = te$ for some idempotent $e$. The subsemigroup $E(S)$ of idempotents of $S$ is commutative, and hence forms a (meet) semilattice for the natural partial order where $e \wedge f := ef$ for $e,f$ in $E(S)$. A helpful introduction to the algebraic theory of inverse semigroups can be found in \cite{LawsonBook98}.

The \textit{left regular representation} $\Lambda : S \ra \el{S}$ of an inverse semigroup $S$ is defined by
$$ \Lambda(a)\delta_b = \left\{\begin{array}{ll}
                        \delta_{ab} & \mbox{if $b \in D_a$} \\
                        0  & \mbox{otherwise}
                               \end{array} \right. $$
for $a,b$ in $S$, where $D_a = \{b : a^* a b = b \} = \{b : a^* a \geq b b^* \}$. We also denote by $\Lambda$ the induced $*$-homomorphism on $C^*(S)$. The algebra $C^*_{r}(S) :=
\Lambda(C^*(S))$ is called the reduced $\cstar$-algebra of $S$. The right regular representation $R$ of $S$ and the induced map on $C^*(S)$ are defined similarly in terms of right multiplication on $\el{S}$. It is easy to show that the images of $\Lambda$ and $R$ commute. We will sometimes write $\Lambda_S$, $R_S$ to avoid confusion when there are multiple inverse semigroups in play.

We regard the $*$-algebra $\bbC S$ both as finitely supported functions $f: S \ra \bbC$ and as finite formal sums over semigroup elements $f = \sum f(s) \, s$. The product of $f,g$ in $\bbC S$ is defined as $(f \cdot g)(a) = \sum_{st = a} f(s)g(t)$ and the involution is defined as $f^* = \sum \overline{f(s)} \, s^*$. It is a dense subalgebra of both the full and reduced $C^*$-algebras of $S$. An element $f \in \bbC S$ is \textit{positive} if $f$ can be expressed as a finite sum of elements of the form $g^*g$, where $g \in \bbC S$. A map between the complex algebras of two inverse semigroups is \textit{positive} if it carries positive elements to positive elements. A \textit{state} on the algebra $\bbC S$ of an inverse semigroup $S$ is a positive linear map $\rho : \bbC S \ra \bbC$, such that
\begin{equation}
\sup \{ {\left|\rho(a)\right|}^2 : a\in \bbC S; \rho(a^* a) \leq 1\}
= 1 \tag{$*$}
\end{equation}
This last condition ensures that states on $\bbC S$ induce cyclic representations via the GNS construction. See section 1 of \cite{DunPat85} for a discussion of states on $\bbC S$, where the authors use a condition different than $(*)$. For the general theory of representable positive linear maps on $*$-algebras, and for the equivalence of the condition in \cite{DunPat85} to $(*)$, see Palmer \cite[Sec. 9.4]{Palmer01}.

If $\rho$ is a state on $\bbC S$, the GNS construction furnishes a $*$-representation $\pi_{\rho} : \bbC S \ra \mcB(\mcH_{\rho})$ with cyclic unit vector $x_{\rho}$. The map $\pi_{\rho}$ extends to a representation $\widetilde{\pi}_{\rho}$ on $C^*(S)$. Then
$\widetilde{\rho}(A) := \inp{\widetilde{\pi}_{\rho}(A)x_{\rho},x_{\rho}}$ defines a state on $C^*(S)$ that extends $\rho$. Conversely, if $\rho$ is a state on $C^*(S)$, then $\rho$ restricts to a state on $\bbC S$. This bijective correspondence between states on $C^*(S)$ and states on $\bbC S$ gives the norm formula
$$\norm{f}_{C^*(S)} = \sup\{\rho(f^* f)^{1/2} : \rho \in \mcS(C^*(S))\} = \sup\{\rho(f^* f)^{1/2} : \rho \in \mcS(\bbC S)\} $$ where $f$ lies in $\bbC S$ and $\mcS(X)$ denotes the set of states for the algebra X.

A \textit{Fell bundle} over a discrete group $G$ is a collection of closed subspaces $\mcB = \{B_g\}_{g \in G}$ of a $C^*$-algebra $B$, satisfying ${B_g}^* = B_{g^{-1}}$ and $B_g B_h \subseteq B_{gh}$ for all $g$ and $h$ in $G$. If the subspaces $B_g$ are linearly independent and their direct sum is dense in $B$, then $\mcB$ is called a grading for $B$. If in addition there is a conditional expectation $\eps : B \ra B_{1_G}$, where $1_G$ is the identity of $G$, then $\mcB$ is called a topological grading for $B$. Exel \cite{Exel97} has defined the reduced $C^*$-algebra $C^*_{r}(\mcB)$ of a Fell bundle $\mcB$, and has shown that all topologically graded $C^*$-algebras over $\mcB$ lie between $C^*(\mcB)$ and the reduced $C^*$-algebra $C^*_{r}(\mcB)$, both of which are graded over $\mcB$. When the two algebras are isomorphic the Fell bundle is said to be \textit{amenable}. This is the case if and only if the expectation on $C^*(\mcB)$ is faithful. For the definition of $C^*(\mcB)$, see \cite[VIII.17.2]{FellDoran88}. We remark that the $*$-representations of $C^*(\mcB)$ are in one-to-one correspondence with the $*$-representations of $\mcB$, and so $C^*(\mcB)$ can be thought of as the full $C^*$-algebra of $\mcB$.

\section{Weak containment and homomorphisms onto groups} \label{secGroupHom}
Let $\varphi : S \ra G$ be a homomorphism of an inverse semigroup $S$ onto a group $G$ with kernel $H$. Under the assumption that $H$ is $C^*$-isometric we show that $S$ has weak containment if and only if $H$ has weak containment and $G$ is amenable. This is the analog of the fact that, in an exact sequence of discrete groups $H \into G \onto K$, $G$ is amenable if and only if $H$ and $K$ are amenable. The main corollary to this fact is that an $E$-unitary inverse semigroup has weak containment if and only if the maximum group image $G(S)$ is amenable.

We first show that the map $\eps : \bbC S \ra \bbC H$ defined by $$ \eps(\sum_{s \in S}{\alpha_s s}) = \sum_{h \in H}{\alpha_h h}$$ extends to a faithful expectation $\eps_r : C^*_{r}(S) \ra \overline{\bbC H}^{{}_{r}}$ and an expectation $\eps_f : C^*(S) \ra \overline{\bbC H}^{{}_{f}}$, where $\overline{\bbC H}^{{}_{r}}$ is the closure of $\bbC H$ in $C^*_{r}(S)$ and $\overline{\bbC H}^{{}_{f}}$ is the closure in $C^*(S)$. The proof is inspired by the work on coactions of groups on $C^*$-algebras appearing in papers such as \cite{LPRS87} and \cite{Quigg96}. The map $\delta$ below is in fact a coaction of $G$ on $C^*_{r}(S)$. The construction of the unitary $W$ is adapted from \cite{LPRS87}, while the proof that $\eps_r$ is faithful mimics the proof of Lemma 1.4 in \cite{Quigg96}.

\begin{proposition}
There exists a faithful conditional expectation $$\eps_r : C^*_{r}(S) \ra \overline{\bbC H}^{{}_{r}}$$ such that $$ \eps_r(\sum_{s \in S}{\alpha_s s}) = \sum_{h \in H}{\alpha_h h}.$$
\end{proposition}

\begin{proof}
When viewing a semigroup element $t$ in $C^*_{r}(S)$ we will now write $\Lambda(t)$. Let $\lambda$ denote the $*$-homomorphism on $C^*(G)$ induced by the left regular representation of $G$. We first show that the map $\Lambda(t) \mapsto \Lambda(t) \otimes \lambda(\varphi(t))$ extends to a $*$-homomorphism $\delta : C^*_{r}(S) \ra C^*_{r}(S) \motimes C^*_{r}(G)$. Since there is no universal property for $C^*_{r}(S)$ this is a nontrivial fact. First, let $W$ in $\mcB( \el{S} \otimes \el{G})$ be the unitary operator defined by $$ W (\delta_{s} \otimes \delta_{g}) = \delta_{s} \otimes \delta_{\varphi(s)g}.$$ Define $\delta := \Ad(W) \circ j$, where $j(A) = A \otimes I$ for $A$ in $C^*_{r}(S)$. Then $\delta$ is a (bounded) $*$-homomorphism. It is also clear that $\delta$ is injective. For $\delta_{s} \otimes \delta_{g}$ in $\el{S} \otimes \el{G}$ and $t$ in $S$ we have:
\begin{eqnarray*}
\delta(\Lambda(t))(\delta_{s} \otimes \delta_{g}) & = & W j(\Lambda(t)) W^* (\delta_{s} \otimes \delta_{g}) \\
                        & = & \left\{\begin{array}{ll}
                                        \delta_{ts} \otimes \delta_{\varphi(t)g} & \mbox{if $s \in D_t$} \\
                                        0  & \mbox{otherwise}
                                     \end{array} \right. \\
                        & = & \Lambda(t) \otimes \lambda(\varphi(t))(\delta_{s} \otimes
                        \delta_{g})
\end{eqnarray*}
Then $\delta$ satisfies $\delta(\Lambda(t)) = \Lambda(t) \otimes \lambda(\varphi(t))$.

Given $\rho \in C^*_{r}(G)^*$ there is a slice map $S_{\rho} : C^*_{r}(S) \motimes C^*_{r}(G) \ra C^*_{r}(S)$. That is, $S_{\rho}$ is a linear map of norm $\norm{\rho}$ such that
$$S_{\rho}( \sum_{i=1}^{n}{a_i \otimes b_i} )= \sum_{i=1}^{n}{a_i \rho(b_i)}$$ for $a_i \otimes b_i$ in $C^*_{r}(S) \motimes C^*_{r}(G)$. Moreover, if $\rho$ is positive, then $S_{\rho}$ is a completely positive map \cite[Corollary IV 4.25]{Tak02}. Recall that $\chi_{e}$, where $e$ is the identity of $G$, is a faithful state on $C_{r}^*(G)$. Thus $ \eps_r = S_{\chi_{e}} \circ \delta $ is a positive contraction on $C^*_{r}(S)$. Also, for $s \in S$,
$$\eps_r(\Lambda(s)) = S_{\chi_e}(\Lambda(s) \otimes \varphi(s)) = \left\{\begin{array}{ll}
   \Lambda(s) & \mbox{if $s \in H$} \\
   0 & \mbox{otherwise}
   \end{array} \right. $$
It follows by linearity that $$\eps_r(\sum_{s \in S}{\alpha_s \Lambda(s)}) = \sum_{h \in H}{\alpha_h \Lambda(h)}$$ and by continuity that $\eps_r$ is a projection with range $\overline{\bbC H}$.

Suppose that $a \neq 0$ in $C^*_{r}(S)$ is positive. Since $\delta$ is injective, $\delta(a)$ is nonzero and positive, and hence there is a state $\omega$ of $C^*_{r}(S)$, such that $(\omega \otimes \iota) \circ \delta(a)$ is a nonzero positive element of $C^*_{r}(G)$, where $\iota$ is the identity on $C^*_{r}(G)$ (c.f. \cite{Wasserman76}). Then $ \chi_{e} \circ (\omega \otimes \iota) \circ \delta(a) = \omega \circ (\iota \otimes \chi_{e}) \circ \delta(a) = \omega \circ \eps_r(a)$. Thus $\eps_r(a) \neq 0$.
\end{proof}

We can construct an expectation $\eps_f$ on $C^*(S)$ in a similar way, using full $C^*$-algebras instead of reduced ones. This map will not be faithful in general since $\chi_{e}$ is not always faithful on $C^*(G)$. However if $G$ is amenable, then $\chi_{e}$ is faithful and hence $\eps_f$ is faithful. These facts are recorded in the next proposition.

\begin{proposition}\label{prop_faithfulexpectation}
There is a conditional expectation $\eps_f : C^*(S) \ra \overline{\bbC H}^{{}_{r}}$ such that $ \eps(\sum_{s \in S}{\alpha_s s}) = \sum_{h \in H}{\alpha_h h}.$ Moreover, $\eps$ is faithful if $G$ is amenable.
\end{proposition}

In order to relate the weak containment property on $S$ to the weak containment property on $H$, we need an assumption about the norms of elements of $\bbC H$ inside $C^*(S)$. Recall the following:

\begin{definition} An inverse subsemigroup $H$ of an inverse semigroup $S$ is \textit{$C^*$-isometric} in $S$ if $\norm{f}_{C^*(S)} = \norm{f}_{C^*(H)}$ for all $f$ in $\bbC H$.
\end{definition}

If $H$ is $C^*$-isometric in $S$ then the embedding of $\bbC H$ in $C^*(S)$ extends to an isomorphism between $C^*(H)$ and $\overline{\bbC H}^{{}_{f}}$. This fact is used implicitly in the proof of the next theorem, the main result of this section.

\begin{theorem}\label{thm_weakcontainment} Let $\varphi : S \ra G$ be a homomorphism of an inverse semigroup $S$ onto a group $G$ with kernel $H$. Suppose that $H$ is $C^*$-isometric in $S$. The following conditions are equivalent:
\begin{itemize}
\item[(i)] $S$ has weak containment,
\item[(ii)] $\eps_f$ is faithful and $H$ has weak containment.
\end{itemize}
\end{theorem}

\begin{proof}
Let $\Lambda_H$ denote the map on $C^*(H)$ induced from the left regular representation of $H$ and let $\Lambda_{S,H}$ denote the restriction of $\Lambda_S$ to $C^*(H)$. We first show that $\Lambda_H$ is injective if and only if $\Lambda_{S,H}$ is injective. For $h \in H$ and $s \in D_h$ notice that, $$ hs \in H \mbox{ if and only if } s \in H.$$ It follows that $\el{H}$ is an invariant subspace for $\Lambda_{S}(h)$, and that
$$\Lambda_{S}(h) =
\begin{bmatrix}
   \Lambda_H(h) & 0 \\
   0            & \; * \;
\end{bmatrix}$$
with respect to the decomposition $\el{H} \oplus \el{H}^{\perp}$. Thus, for any $A$ in $C^*(H)$,
$$\Lambda_{S}(A) =
\begin{bmatrix}
   \Lambda_H(A) & 0 \\
   0            & \; * \;
\end{bmatrix}$$
Hence if $\Lambda_H$ is injective, then $\Lambda_{S,H}$ is injective. Conversely, if $\Lambda_{S}(A) \neq 0$ for some $A$ in $C^*(H)$, then $\Lambda_{S}(A) \, \delta_{s} \neq 0$ for some $s$ in $S$. But then,
$$ R_{S}(s) \Lambda_{S}(A) \, \delta_{s s^*} = \Lambda_{S}(A)
R_{S}(s) \, \delta_{s s^*} = \Lambda_{S}(A) \, \delta_{s} \neq 0$$
Since $s^* s \in H$ we have that $\Lambda_{H}(A) \neq 0$.

We see from the commuting diagram:
\[
\begin{CD}
C^*(S)   @>\Lambda_{S}>>    C^*_{r}(S) \\
@V\eps_f VV                        @VV\underset{\text{faithful}}{\eps_{r}}V \\
C^*(H)   @>\Lambda_{S,H}>>  \overline{\bbC H}
\end{CD}
\]
that $\Lambda_{S}$ is injective if and only if $\eps_f$ is faithful and $\Lambda_{S,H}$ is injective. Since $\Lambda_{S,H}$ is injective if and only if $\Lambda_{H}$ is injective, the theorem follows.
\end{proof}

Paterson has shown that $G(S)$ is amenable when $S$ has weak containment (See \cite[Proposition 4.1]{Pat78} and note that, in Paterson's notation, $1 \in P(S) = P_{L}(S)$). Since $G(S)$ maps onto $G$, it follows that $G$ is amenable if $S$ has weak
containment. Also, by Proposition \ref{prop_faithfulexpectation}, $\eps$ is faithful if $G$ is amenable. Combining these two facts with the theorem we get the following corollary.

\begin{corollary}\label{cor_weakcontainment}
Let $\varphi : S \ra G$ be a homomorphism of an inverse semigroup $S$ onto a group $G$ with kernel $H$. Suppose that $H$ is $C^*$-isometric in $S$. The following conditions are equivalent:
\begin{itemize}
\item[(i)] $S$ has weak containment,
\item[(ii)] $G$ is amenable and $H$ has weak containment.
\end{itemize}
\end{corollary}

Though it is often the case that $H$ is $C^*$-isometric, we do not know if this property always holds. One way to show $\norm{f}_{C^*(S)} = \norm{f}_{C^*(H)}$ for all $f$ in $\bbC H$ is to prove $\eps : \bbC S \ra \bbC H$ is positive. Suppose $\eps$ is positive. If $\rho$ is a state on $\bbC H$, then $\rho \circ \eps$ is a positive linear map. In fact, for $a$ in $\bbC S$, $\eps(a)^* \eps(a) \leq \eps( a^* a)$. Then for each $a$ in $\bbC S$ with $(\rho \circ \eps)(a^*a) \leq 1$, $\eps(a) \in \bbC H$ with $\rho(\eps(a)^* \eps(a)) \leq 1$. Thus $\rho \circ \eps$ satisfies condition ($*$) from section 1, and is therefore a state on $\bbC S$. Fix $f$ in $\bbC H$. Any representation of $S$ restricts to a representation of $H$, from which it follows that $\norm{f}_{C^*(H)} \geq \norm{f}_{C^*(S)}$. Conversely, since states on $\bbC H$ can be extended to states on $\bbC S$,
\begin{eqnarray*}
 \norm{f}_{C^*(H)} & = & \sup\{\rho(f^* f)^{1/2} : \rho \in \mcS(\bbC H)\} \\
                      & \leq & \sup\{\rho(f^* f)^{1/2} : \rho \in \mcS(\bbC S)\} \\
                      & = & \norm{f}_{C^*(S)}
\end{eqnarray*}
Thus, $\norm{f}_{C^*(H)} = \norm{f}_{C^*(S)}$.

In section 6, we show there are many inverse semigroups $S$ for which $\eps$ is positive. Most notably, if $S$ is $E$-unitary then the kernel of the homomorphism $\sigma : S \ra G(S)$ is the semilattice $E := E(S)$. We show $\eps : \bbC S \ra \bbC E$ is positive. Since any semilattice has weak containment we get the following result.

\begin{corollary}\label{cor_eunitary} Suppose $S$ is an $E$-unitary inverse semigroup. Then $S$ has weak containment if and only if $G(S)$ is amenable.
\end{corollary}

\section{Weak containment and inverse semigroups with zero} \label{secGroupGrading}

The use of group homomorphisms in the previous section is suitable for many inverse semigroups, including all $E$-unitary inverse semigroups. However, if $S$ contains a zero, then $G(S)$ is trivial, $S$ is the kernel of any homomorphism onto a group, and the results of the previous section are vacuous. To remedy this, we work with maps that are not quite homomorphisms. For any inverse semigroup $S$, let $S^0$ denote the inverse semigroup obtained from $S$ by adjoining a zero if $S$ does not already have one, otherwise $S^0 = S$.

\begin{definition}
A \textit{grading} of an inverse semigroup $S$ containing a zero by the group $G$ is a map $\varphi : S \ra G^0$ such that $\varphi^{-1}(0) = \{0\}$ and $\varphi(ab) = \varphi(a) \varphi(b)$ provided that $ab \neq 0$.
\end{definition}

It is customary when working with algebras generated by semigroups with zero to consider the quotient by the ideal generated by the zero. This identifies the zero of the algebra with the zero of the semigroup. The algebras $C^*_0(S), C^*_{r \, 0}(S)$, and $\bbC_0 S$, for example, are just the quotients of the algebras with which we have been working by the ideal generated by the zero of $S$.

Fix a grading $\varphi : S \ra G^0$, and let $H = \varphi^{-1}(1_G)^0$. As before, we have an expectation $\eps_f : C^*_0(S) \ra \overline{\bbC_{0} H}^{{}_{f}}$ and a faithful expectation $\eps_r : C^*_{r \, 0}(S) \ra \overline{\bbC_{0} H}^{{}_{r}}$, where each map extends the restriction map from $\bbC_0 S$ onto $\bbC_0 H$. The analog of Theorem \ref{thm_weakcontainment} holds.

\begin{theorem}\label{thm_weakcontainmentzero} Let $\varphi : S \ra G^0$ be a grading of an inverse semigroup $S$ by a group $G$. Let $H = \varphi^{-1}(1_G)^0$. Suppose that $H$ is $C^*$-isometric in $S$. Then $S$ has weak containment if and only if $\eps_f$ is faithful and $H$ has weak containment.
\end{theorem}

Unlike in the previous section, it may happen that $S$ has weak containment and yet $G$ is not amenable. In the next section we show that all graph inverse semigroups have weak containment, yet the universal grading of a graph inverse semigroup is a free group. The proof requires that we view a grading $\varphi : S \ra G^0$ in a different light. We observe that $\varphi$ induces a Fell bundle structure (see section 1) on $C^*_0(S)$. We then relate weak containment for $S$ to amenability of the Fell bundle, showing in particular that strongly $E^*$-unitary inverse semigroups (such as graph inverse semigroups) have weak containment if and only if the associated Fell bundle is amenable. Exel \cite{Exel97} has found an approximation property for Fell bundles that guarantees amenability. Using this property he was able to give examples of amenable Fell bundles over nonamenable groups arising from Cuntz-Krieger algebras. We use conditions on Fell bundles over free groups that Exel found in subsequent work \cite{Exel00} to establish amenability of Fell bundles arising from graph inverse semigroups.

We first define the Fell bundle structure arising from a grading $\varphi$. For each $g$ in $G$, let
\begin{eqnarray*}
A_g & = & \Span\{s : \varphi(s) = g\} \mbox{ inside } {\bbC_0 S} \\
B_g & = & \overline{A_g} \mbox{  inside } C^*_0(S) \\
\end{eqnarray*}

\begin{proposition}
The collection $\mcB = \{B_g\}_{g \in G}$ is a Fell bundle for $C^*_0(S)$.
\end{proposition}
\begin{proof}
We show only that $A_g A_h \subseteq A_{gh}$ for all $g$ and $h$ in $G$. It is enough to show that, for $s$ and $t$ in $S$ such that $\varphi(s) = g$ and $\varphi(t) = h$, $st \in A_{gh}$. This is the case since either $st = 0 \in A_{gh}$, or $\varphi(st) = \varphi(s) \varphi(t) = gh$, in which case $st \in A_{gh}$.
\end{proof}

In fact, since there is an expectation $\eps_f : C^{*}_0(S) \ra B_{1_G}$ that vanishes on $B_{g}$ for $g \neq 1_G$, $\mcB$ is a topological grading for $C^{*}_0(S)$ \cite[Theorem 3.3]{Exel97}. It follows that $\eps_f$ is faithful whenever $\mcB$ is amenable. In fact, one can verify that representations of $\mcB$ are in one-to-one correspondence with representations of $C^{*}_0(S)$ and hence $C^{*}(\mcB)$ is isomorphic to $C^{*}_0(S)$. Moreover, the expectation on $C^{*}(\mcB)$ is just $\eps_f$. Thus $\eps_f$ is faithful if and only if $\mcB$ is amenable. Recall that a strongly $E^{*}$-unitary inverse semigroup $S$ is an inverse semigroup that admits a grading $\varphi : S \ra G^{0}$, with
$\varphi^{-1}(1_G)$ equal to the nonzero idempotents of $S$. We then have the following corollary to Theorem \ref{thm_weakcontainmentzero}, which is the analog of Corollary \ref{cor_eunitary} for inverse semigroups with zero.

\begin{corollary} \label{cor_strongeunitary}Suppose $S$ is strongly $E^{*}$-unitary. Then
$S$ has weak containment if and only if the Fell bundle $\mcB$ arising from the universal grading of $S$ is amenable.
\end{corollary}

\section{Graph inverse semigroups} \label{secGraphSemigroup}

We follow most of the conventions of \cite{Rae05} for directed graphs. Briefly, a \textit{directed graph} $\msE = (\msE^0, \msE^1, r, s)$ consists of countable sets $\msE^0$, $\msE^1$ and functions $r,s : \msE^1 \ra \msE^0$. The elements of $\msE^0$ are called \textit{vertices}, and the elements of $\msE^1$ are called \textit{edges}. Given an edge $e$, $r(e)$ denotes the range of $e$ and $s(e)$ denotes the source of $e$. We denote by $\msE^*$ the collection of finite directed paths in $\msE$. The functions $r,s$
can be extended to $\msE^*$ by defining $r(\mu) = r(\mu_n)$, $s(\mu) = s(\mu_1)$ for a path $\mu = \mu_n \mu_{n-1} \cdots \mu_1$ in $\msE^*$. If $\mu = \mu_n \mu_{n-1} \cdots \mu_1$ and $\nu = \nu_m \nu_{m-1} \cdots \nu_1$ are paths with $s(\mu) = r(\nu)$, we write
$\mu \nu$ for the path $\mu_n \cdots \mu_1 \nu_m \cdots \nu_1$. The length of a path $\mu$ is denoted $\left|\mu\right|$.

The \textit{graph inverse semigroup} of the directed graph $\msE$ is the set
$$ S_\msE = \{(\mu,\nu) : s(\mu) = s(\nu) \} \cup \{ 0 \} $$
with the products not involving zero defined by
$$ (\mu,\nu) (\alpha,\beta) = \left\{\begin{array}{ll}
        (\mu, \beta \nu') & \mbox{if $\nu = \alpha \nu'$} \\
        (\mu \alpha',\beta) & \mbox{if $\alpha = \nu \alpha'$} \\
        0 & \mbox{otherwise}
   \end{array} \right. $$
The inverse operation is given by $(\mu,\nu)^{*} = (\nu,\mu)$. It is easy to see that the set of idempotents of $S_{\msE}$ is $E = \{(\mu,\mu): \mu \in \msE^* \}$.

The inverse semigroup $S_{\msE}$ is important in the study of $C^*$-algebras of directed graphs. It has been shown that $C^*(\msE)$ is a quotient of $C^*(S_{\msE})$ \cite{Pat02}. $S_{\msE}$ has also been studied in the semigroup literature. See \cite{Lawson01}, \cite{AshHall75}, for example.

Let $\bbF$ be the free group generated by the set $\msE^1$. Define a map $\varphi : S_{\msE} \ra \bbF$ by $\varphi( (\mu,\nu ) ) = \red(\mu \nu^{-1})$, where $\red(w)$ denotes the reduction of the word $w$ over the alphabet  $\msE^1 \cup (\msE^1)^{-1}$ in $\bbF$. Then $\varphi$ is a grading of $S_{\msE}$ by $\bbF$ with kernel $E$ \cite{Lawson01}. Let $\mcB := \{B_w\}_{w \in \bbF}$ be the Fell bundle for $C^*_{0}(S_{\msE})$ arising from the grading $\varphi$. We want to show that $S_{\msE}$ has weak containment. Since $S_{\msE}$ is strongly $E^*$-unitary, it suffices to show that $\eps: C^*_{0}(S_{\msE}) \ra C^*_{0}(E)$ is faithful. That is, it suffices to show that $\mcB$ is amenable. We first need some definitions.

\begin{definition}[Exel] A Fell bundle $\mcB = \{B_w\}_{w \in \bbF}$ over a free group $\bbF$ with a fixed set of generators $X$ is \textit{orthogonal} if $B_{x}^* B_y = 0$ for distinct $x,y \in X$. $\mcB$ is \textit{semi-saturated} if, for any pair $s,t$ in $\bbF$
such that the product $st^{-1}$ does not involve cancellation, $B_{st^{-1}} = B_s B_{t^{-1}}$.
\end{definition}

The following theorem was proved in \cite{Exel00}.
\begin{theorem}[Exel]
Let $\mcB$ be an orthogonal, semi-saturated Fell bundle over a free group $\bbF$ with separable fibers. Then $\mcB$ is amenable.
\end{theorem}

The rest of this section is devoted to proving the following theorem by showing that the Fell bundle $\mcB$ arising from a graph inverse semigroup is orthogonal and semi-saturated.
\begin{theorem} The inverse semigroup $S_{\msE}$ of a directed graph $\msE$ has weak containment.
\end{theorem}
\begin{proof}
For $w$ in $\bbF$ let $S_w = \varphi^{-1}(w)$. Take $x,y$ in $\msE^1$. An arbitrary element of $S_{x}^*$ is of the form $(\mu, x \mu)$ with $\mu$ in $\msE^*$. Similarly, the elements of $S_{y}$ are of the form $(y \nu,\nu)$, where $\nu$ in $\msE^*$. The product $(\mu, x \mu)(y \nu,\nu)$ = 0 unless either $x \mu = y \nu \alpha$, or $y \nu = x \mu \alpha$ for some $\alpha$ in $\msE^*$. In either case, $x=y$. It follows that $\mcB$ is orthogonal.

Next we show that $\mcB$ is semi-saturated. Suppose $s,t \in \bbF$ where the product $st^{-1}$ involves no cancellation. If $B_{st^{-1}}$ is the zero subspace then the containment $B_{s} B_{t^{-1}} \subseteq B_{st^{-1}}$ implies the two subspaces are equal. Otherwise $st^{-1} = \red(ab^{-1})$ where $a,b$ describe paths in $\msE$ starting at a common vertex $v$. We may assume that $ab^{-1}$ is a reduced word. Also, since $a,b$ are positive words over the set $\msE^1$, either $a$ is a prefix of $s$ or $b$ is a prefix of $t$. Since the two cases are similar we consider only the first. We can then write $s = a c^{-1}$, where $c$ is a path starting at $v$. Let $E^{v}$ denote the set of idempotents $(w,w)$ with $r(w) = v$. We claim that
$$S_{st^{-1}} = (a,v)E^{v}(v,b).$$
The nontrivial inclusion is $S_{st^{-1}} \subseteq (a,v)E^{v}(v,b)$. Suppose $\varphi((\alpha,\beta)) = st^{-1}$. Then $\red(\alpha \beta^{-1}) = a b^{-1}$. Since $\alpha, \beta$ are already reduced words, the only cancellation in the product $\alpha \beta^{-1}$ occurs where $\alpha$ meets $\beta^{-1}$. Hence, there exists a path $w$ with $r(w) = v$ such that $\alpha = aw$ and $\beta = bw$. Thus $(\alpha, \beta) = (a,v)(w,w)(v,b) \in (a,v)E^{v}(v,b).$

An element $f$ in the span of $S_{st^{-1}}$ can be written
$$ f = \sum_{w \in E^{v}} \lambda_w \, (a,v) (w,w) (v,b), $$
where all but finitely many of the $\lambda_w$ are zero. Suppose $w,w'$ are paths of the same length with $r(w) = v = r(w')$. The product $(w,w)(w',w')$ is nonzero only if $w=w'$, in which case it is $(w,w)$. Set $E^{v}_{k} := \{w : r(w)=v, \left|w\right| = k\}$
and define
\begin{eqnarray*}
f_k   & := & \sum_{w \in E^{v}_{k}} \lambda^{\frac{1}{2}}_w \, (a,v)(w,w)(v,c)\\
f'_k  & := & \sum_{w' \in E^{v}_{k}} \lambda^{\frac{1}{2}}_{w'} \,(c,v)(w',w')(v,b)
\end{eqnarray*}
For $w,w'$ in $E^{v}_{k}$ we have
$$ \left[(a,v)(w,w)(v,c)\right]\left[(c,v)(w',w')(v,b)\right]  = \left\{\begin{array}{ll}
                        (a,v)(w,w)(v,b) & \mbox{if $w = w'$} \\
                        0  & \mbox{otherwise}
                        \end{array} \right. $$
We then have
$$ f = \sum_{k} f_k f'_k.$$

Since $f_k \in B_{s}$ and $f'_k \in B_{t^{-1}}$, we have $f$ in $B_s B_{t^{-1}}$ and it follows that $\mcB$ is semi-saturated.
\end{proof}

\section{Nica's inverse semigroup} \label{secNicaSemigroup}
In \cite{Nica92}, Nica studies $C^*$-algebras $C^*(G,P)$ and $\mcW(G,P)$ associated with certain pairs $(G,P)$ called quasi-lattice ordered groups. Here $G$ is a discrete group and $P$ is a subsemigroup of $G$. He defines $(G,P)$ to be \textit{amenable}
if and only if a natural map $C^*(G,P) \ra \mcW(G,P)$ is an isomorphism. The $C^*$-algebras constructed from quasi-lattice ordered groups include many famous $C^*$-algebras having certain uniqueness properties. It is shown that the uniqueness property follows from amenability of the quasi-lattice ordered group. The first example is $(\bbZ,\bbN)$, from which one recovers the $C^*$-algebra of the unilateral shift. The pair $(\bbZ,\bbN)$ is amenable and this corresponds to the uniqueness property given by
Coburn's theorem. In a subsequent paper \cite{Nica94}, Nica studies an inverse semigroup $\mcT_{G,P}$ induced from a quasi-lattice ordered group $(G,P)$. For example, $\mcT_{\bbZ,\bbN}$ is isomorphic to the bicyclic monoid, an inverse semigroup that has the weak containment property.

In this section, we point out that Nica's definition of amenability of a quasi-lattice ordered group $(G,P)$ is equivalent to weak containment for $\mcT_{G,P}$. For the rest of this section we consider a pair $(G,P)$, where $G$ is a discrete group with a subsemigroup $P$, such that $P \cap P^{-1}$ is the unit of $G$. It follows that the relation $\leq$, defined by $x \leq y$ if and only if $x^{-1}y \in P$, is a partial order on $G$.

\begin{definition}[Nica]
The group $(G,P)$ is \textit{quasi-lattice ordered} if and only if $(1)$ Any $x \in PP^{-1}$ has a least upper bound in P, and $(2)$ Any $s,t \in P$ with a common upper bound have a least common upper bound.
\end{definition}

One often considers the inverse semigroup $I(X)$ of partially-defined bijections on the set $X$. That is, a function $f$ in $I(X)$ is a bijection of a subset $\dom(f)$ of $X$ to another subset $\ran(f)$ of $X$. The multiplication of elements $f,g$ in $I(X)$ is given by composition of the two functions on the largest domain where the composition is defined. The inverse semigroup $I(X)$ plays the same role in inverse semigroup theory as the group of permutations on a set plays in group theory (c.f. \cite[Theorem 1, p. 36]{LawsonBook98}).

For each $x \in G$, define $\beta_{x} : \{t \in P : xt \in P\} \ra \{s \in P : x^{-1}s \in P\}$ by $\beta_{x}t = xt$. Notice that $\beta_{x} \in I(P)$. Then $\mcT_{G,P}$ is defined to be the inverse semigroup generated by ${\{\beta_x\}}_{x \in G}$. It is shown in \cite[Theorem 6.9]{Lawson01}, that $\mcT_{G,P}$ is strongly $E^*$-unitary. An idempotent pure grading $\varphi : \mcT_{G,P} \ra G$ is given by $\varphi(\beta_{x_1} \dots \beta_{x_n}) = x_1 \dots x_n$ for $\beta_{x_1} \dots \beta_{x_n} \neq 0$. Thus, by Corollary \ref{cor_strongeunitary}, $\mcT_{G,P}$ has weak containment if and only if the conditional expectation $\eps : C^{*}_{0}(\mcT_{G,P}) \ra C^{*}_{0}(E)$ is faithful, where $E = E(\mcT_{G,P})$.

Notice that the semigroup isomorphism given near the end of page $370$ in \cite{Nica94} shows that the algebra $C^{*}(G,P)$ defined in \cite[Section 4.1]{Nica92} is isomorphic to
$C^{*}_{0}(\mcT_{G,P})$. Moreover, the conditional expectation on $C^{*}(G,P)$ is $\eps$. Thus, by the first proposition in Section 4.3 of \cite{Nica92}, we have:

\begin{proposition} A quasi-lattice ordered group $(G,P)$ is amenable if and only if the inverse semigroup $\mcT_{G,P}$ has weak containment.
\end{proposition}

\section{Positivity of $\eps$} \label{secPositivityofEpsilon}

Let $H$ be an inverse subsemigroup of an inverse semigroup $S$, and $\eps : \bbC S \ra \bbC H$ the restriction map. The crucial hypothesis that $H$ is $C^*$-isometric in $S$ that appears in section 2 is satisfied when $\eps$ is positive. In this section we study that positivity property. It is easy to find examples where $\eps$ fails to be positive for an arbitrary subsemigroup $H$. We give one such example below (Example \ref{exEpsNotPos}). We also find two classes of semigroups where it is possible to prove positivity. The question raised in this section, which we do not answer, is whether $\eps$ is positive when $H$ is the kernel of a group homomorphism. There is some reason to believe the question has an affirmative answer. In the case that $S$ is a group, the restriction map onto the complex algebra generated by any subgroup is positive. This was proved by Rieffel \cite[Lemma 1.1]{Rie74} using basic facts about cosets of $H$ in $S$. In general, an inverse subsemigroup $H$ does not admit cosets that partition the larger semigroup. However, if $H$ is closed upwards in the natural partial order (such as the case that $H$ is the kernel of a group homomorphism), then there is a related notion of \textit{$\omega$-cosets}, defined by
$$ \uparrow \! sH := \{t \in S : te \in sH, \mbox{ for some idempotent } e\} $$
for $s$ in $S$. These sets partition $S$ and play a role similar to that of cosets of a subgroup (c.f. \cite[IV.4]{PetrichBook}). We suspect it may be possible to generalize Rieffel's proof to such subsemigroups, but so far we have not been successful. For this reason we pose the following question, which is more general than the case that $H$ is the kernel of a homomorphism onto a group.

\begin{question} Suppose $H$ is a inverse subsemigroup of $S$ such that $\ua H = H$. Does it follow that the restriction map $\eps : \bbC S \ra \bbC H$ is positive?
\end{question}

We now give an example showing that $\eps$ is not positive if the hypothesis that $\ua H = H$ is not satisfied.

\begin{example} \label{exEpsNotPos}
Let $a$ in $I(\bbZ)$ (see section 5) be the map defined by $a(n) = n + 1$ for $n$ in $\bbZ$, and let $e$ in $I(\bbZ)$ be the identity on $\bbN$ and undefined elsewhere. Let $S$ be the inverse semigroup generated by $a$ and $e$ and let $H$ be the inverse subsemigroup of $S$ generated by $b := ae$. Then $H$ is the bicyclic semigroup and $C^*(H)$ is the $C^*$-algebra of the unilateral shift. Let $x = e - a$. Then $x x^*$ is a positive element of $\bbC S$ and $ \eps(x x^*) = e - b - b^*$. This is not a positive element of $\bbC H$. To see this note that $e - b - b^*$ maps to the non-positive function $1 - z - \overline{z}$ in the Calkin algebra $C(\bbT)$.
\end{example}

We can prove positivity of $\eps$ if either $H$ is the semilattice of $S$, or if, for each $\omega$-coset $\ua s H$, there is an element $s'$ such that $\ua s H = s' H$. The next result will be used or that purpose.

\begin{proposition} Let $\varphi : S \ra G$ be a homomorphism of an inverse semigroup $S$ onto a group $G$ and let $H = \varphi^{-1}(1_G)$. Suppose $\eps : \bbC S \ra \bbC H$ is the restriction map. Then $\eps(f^*f)$ is positive in $\bbC S$ for every $f$ in $\bbC S$.
\end{proposition}

\begin{proof} We can write
$$ f = \sum_{g \in G} f_g $$
where $f_g \in \Span{\varphi^{-1}(g)}$. Then
$$\eps(f^*f) = \eps(\sum_{g,h \in G} {f_h}^* f_g) = \sum_{g \in G} {f_g}^* f_g, $$
a positive element of $\bbC S$.
\end{proof}

\begin{remark}
Note we have not shown that $\eps$ is positive, since we do not know that $\eps(f^* f)$ is positive as an element of $\bbC H$. The author would like to thank an anonymous referee for finding this error in a previous version of this paper.
\end{remark}

\begin{proposition} \label{eps_positive} Let $\varphi : S \ra G$ be a homomorphism of an inverse semigroup $S$ onto a group $G$ and let $H = \varphi^{-1}(1_G)$. Suppose $\eps : \bbC S \ra \bbC H$ is the restriction map. If $H = E(S)$, or if for each $g \in G$ there is an element $s_g \in S$ such that $\varphi^{-1}(g) = s_g H$, then $\eps$ is positive.
\end{proposition}

\begin{proof}
For $g \in G$ let $S_g := \varphi^{-1}(g)$. By the proof of the last proposition it is enough to show, for every $f \in \Span \{S_g\}$ there exists $f' \in \bbC H$ such that ${f'}^* f' = f^* f$. To this end, fix $f = \sum_{s \in S_g} \alpha_s s$ in $\Span\{S_g\}$. In the case that $H = E(S)$ let
$$ f' = \sum_{s \in S_g} \alpha_s s^* s.$$
Then $f' \in \bbC H$ and, since $s^*(st^*)t = s^*(ts^*)t = s^* t$ for all $s,t \in S_g$ we have:
$$ {f'}^* f' = \sum_{s,t \in S_g} \overline{\alpha_s}\alpha_t s^*st^*t = \sum_{s,t \in S_g} \overline{\alpha_s}\alpha_t s^* t = f^* f.$$

If we have $S_g = s_g H$ for some $s_g \in S$ then let $f' = \sum_{s \in S_g} \alpha_s s_g^* s_g h_s$, where $h_s$ is chosen so that $s = s_g h_s$. Since $(s_g^* s_g)^2 = s_g^* s_g$ we have
$$ {f'}^* f' = \sum_{s,t \in S_g} \overline{\alpha_s}\alpha_t h_s^* (s_g^* {s_g})^2 h_t  = \sum_{s,t \in S_g} \overline{\alpha_s}\alpha_t s^* t = f^* f.$$
\end{proof}

Next we give an example satisfying the second hypothesis of the above proposition that is in general not $E$-unitary. Using the main result of section 2, we can characterize the weak containment property for the class of \textit{bisimple inverse $\omega$-semigroups}. This is the class of bisimple inverse semigroups whose idempotents form a descending chain $e_0 > e_1 > e_2 > \dots$.

\begin{example} Let $\theta$ be an endomorphism of a group $G$, and $BR(G,\theta)$ be the set $\bbN^{0} \times G \times \bbN^{0}$ equipped with the multiplication:
$$ (m, a, n)(i, b, j) = (m - n + t, \theta^{t-n}(a)\theta^{t-i}(b), j - i + t),$$ where $t = \max(n,i)$. $BR(G,\theta)$ is called the \textit{Bruck-Reilly extension of $G$ determined by $\theta$}. Reilly \cite{Reilly66} showed that $S$ is a bisimple inverse $\omega$-semigroup if and only if $S$ is isomorphic to some $BR(G,\theta)$. There is a homomorphism $\varphi: BR(G,\theta) \ra \bbZ$ given by $\varphi(m,a,n) = m - n$. Let $H = \ker \varphi$. For $k \in \bbZ$, notice that
$$\varphi^{-1}(k) = \left\{\begin{array}{ll}
                        (k,1_G,0)H & \mbox{if $k \geq 0$} \\
                        (0,1_G,k)H & \mbox{otherwise}
                    \end{array} \right.$$
Thus, by Proposition \ref{eps_positive}, $H$ is $C^*$-isometric in $BR(G,\theta)$. It follows by Corollary \ref{cor_weakcontainment} that $BR(G,\theta)$ has weak containment if and only if $H$ has weak containment. Since the idempotents of $H$ are central, $H$ is a Clifford $\omega$-semigroup with all maximal subgroups isomorphic to $G$. By \cite[Theorem 2.6]{DunPat85} $H$ has weak containment if and only if $G$ is amenable. Thus
$BR(G,\theta)$ has weak containment if and only if $G$ is amenable. It should be noted that a stronger result is proved by Duncan and Paterson in \cite{DunPat85}, where the authors characterize weak containment for the Bruck-Reilly extension of a finite semilattice of groups.
\end{example}

Finally, we note that there is a version of Proposition \ref{eps_positive} in the case that $S$ contains a zero. We state it here without proof.

\begin{proposition} Let $\varphi : S \ra G^0$ be a grading of an inverse semigroup $S$ containing a zero by a group $G$ and let $H = \varphi^{-1}(1_G)^0$. Suppose $\eps : \bbC_0 S \ra \bbC_0 H$ is the restriction map. If $H = E(S)$, or if for each $g$ in $G$ there is an element $s_g$ in $G$ such that $\varphi^{-1}(g) = s_g H$, then $\eps$ is positive.
\end{proposition}

\begin{acknowledge} This work was completed as part of the author's doctoral dissertation. He gratefully acknowledges the abundant support of his department and his co-advisors, John Meakin and David Pitts.
\end{acknowledge}

\end{document}